\documentstyle[12pt]{article}
  
\begin{titlepage}

\title{The $\eta$-form and a generalized Maslov index}

\author{Ulrich Bunke\thanks{Mathematisches Institut, Universit\"at G\"ottingen, Bunsenstr. 3-5, 37073 G\"ottingen, GERMANY, E-mail:bunke@cfgauss.uni-math.gwdg.de} and H. Koch\thanks{Institut f\"ur Angewandte Mathematik, 
Universit\"at Heidelberg,
Im Neuenheimer Feld 294,
69120 Heidelberg,
GERMANY, E-mail:koch@IWR.Uni-Heidelberg.De
}}

\end{titlepage}

\newcommand{\proof}{{\it Proof.$\:\:\:\:$}}

\newcommand{\R}{{\bf R}}

\newcommand{\Z}{{\bf Z}}
\newcommand{\C}{{\bf C}}

\newcommand{\ch}{{\bf ch}}

\newcommand{\res}{{\rm res}}
\newcommand{\Tr}{{\rm Tr}}

\newcommand{\cD}{{\cal D}}
\newcommand{\cH}{{\cal H}}

\newcommand{\Hom}{{\mbox{\rm Hom}}}

\newcommand{\End}{{\mbox{\rm End}}}

\newcommand{\Ree}{{\rm Re }}

\newcommand{\ee}{{\rm e}}

\newcommand{\tr}{{\mbox{\rm tr}}}

\newcommand{\coker}{{\rm coker}}

\def\hB{\hspace*{\fill}$\Box$ \newline\noindent}

\newcommand{\ind}{{\rm index}}

\def\hB{\hspace*{\fill}$\Box$ \\[0.5cm]\noindent}
\newcommand{\cL}{{\cal L}}

\newtheorem{prop}{Proposition}[section]
\newtheorem{lem}[prop]{Lemma}
\newtheorem{ddd}[prop]{Definition}
\newtheorem{theorem}[prop]{Theorem}

\begin{document}

\maketitle

\tableofcontents
       
\section{Introduction}
 
\newcommand{\dom}{{\rm dom}}
\newcommand{\dirac}{{\not\partial}}
\newcommand{\Aut}{{\rm Aut}}
In this note we consider 
the $\eta$-form of a 
family of Dirac operators $\cD(b)$, $b\in B$,
on the interval $[0,1]$ over a base space
$B$.  The $\eta$-form was introduced by Bismut-Cheeger \cite{bismutcheeger90}
as the boundary contribution to the local index theorem for families of Dirac operators.
In our case the operator $\cD(b)$ depends on $b\in B$ only through
the boundary conditions. 
If $B$ is a point, 
then the $\eta$-form reduces to the usual $\eta$-invariant of $\cD$ which
was explicitly calculated by
Lesch-Wojciechowski \cite{leschwojciechowski96}.
In \cite{bunke95} we found a relation between the $\eta$-invariant and the
Maslov index. The Maslov index was first introduced in Wall \cite{wall69}
as a measure of the non-additivity under gluing of the signature of manifolds with boundary. This non-additivity was generalized to arbitrary
Dirac operators in \cite{bunke95}.

In the present note we relate the $\eta$-form with a family version of the Maslov index. The family version of the Maslov index conjecturally plays the
same role in the non-additivity of the family index of families of Dirac operators on manifolds with cylindrical ends (or APS-boundary conditions) as the usual Maslov index does for the usual index.

Let $V$ be a finite-dimensional Hilbert space equipped with a 
hermitean symplectic structure $\Omega$
(this just means that $\imath\Omega$ is a non-degenerate hermitean form
of index $(l,l)$, $\dim_\C(V)=2l$).
If $\{L_0(b),L_1(b)\}_{b\in B}$
is a smooth family of pairs of transverse Lagrangian subspaces
of $V$, then we define the $\eta$-form $\eta(L_0,L_1)\in C^\infty(B,\Lambda^{ev}T^ *B)$.
Our main result is
\begin{theorem}
If $\{L_0(b),L_1(b),L_2(b)\}_{b\in B}$ is a smooth family 
of triples of pairwise transverse Lagrangian subspaces of $V$,
then\\ 
{\bf (1)}$\:\:$ $d(\eta(L_0,L_1)+\eta(L_1,L_2)+\eta(L_2,L_0))=0$,\\
{\bf (2)}$\:\:$ and if we define the cohomology class $\tau(L_0,L_1,L_2)$ by
 $$\tau(L_0,L_1,L_2):=[\eta(L_0,L_1)+\eta(L_1,L_2)+\eta(L_2,L_0)]\in H^{even}(B,\R)\ ,$$
then $\tau(L_0,L_1,L_2)=\ch(L_0^+)-\ch(L_0^-)$, where
$L_0=L_0^+\oplus L_0^-$ is the splitting of the bundle of Lagrangian subspaces
$L_0\subset B\times V$ into the positive and negative eigenspaces of the quadratic form $Q(x_0):=\Omega(x_1,x_2)$, where $x_i\in L_i$, $x_0=x_1+x_2$.    
\end{theorem}
The proof of the theorem is based on a local index theorem
for families of Dirac operators on manifolds with cylindrical ends
and boundaries with local boundary condition.
Instead of saying how the existing proofs Bismut-Cheeger \cite{bismutcheeger90}, Melrose-Piazza \cite{melrosepiazza96} should be modified
in order to include the additional boundaries
we prefer to work out again the essential arguments.
Our approach is modelled on the $b$-calculus proof
of \cite{melrosepiazza96}, but is more direct and might be of independent interest.

\section{Definition of $\eta(L_0,L_1)$}

Let $V$ be a finite-dimensional complex Hilbert space
with scalar product $(.,.)$.
Let $I\in \End(V)$  be a complex structure, i.e. $I^2=-1$, $I^*=-I$.
We assume that $\tr\:I=0$. Then $\Omega(x,y):=(Ix,y)$ is a hermitean symplectic structure on $V$.
We consider the formally selfadjoint differential operator 
$D:=I\frac{d}{d t}$
acting on $C^\infty([0,1],V)$. 

A complex subspace $L\subset V$ is called Lagrangian, if $L\perp IL$ and $L\oplus IL=V$.
We want to consider $D$ as an unbounded operator on the Hilbert space $L^2([0,1],V)$.
In order to define selfadjoint extensions $\cD$
of $D$ we choose two Lagrangian subspaces $L_i\subset V$, $i=0,1$.
Then we define the domain $\dom(\cD)$ of $\cD$ by
$$\dom(\cD):=\{f\in C^\infty([0,1],V)\:|\:f(i)\in L_i,\:\:i=0,1\}\ .$$  
It can be shown that $\cD$ is essentially selfadjoint and we denote its
unique selfadjoint extension by $\cD$, too.
It is easy to see that $\cD$ is invertible iff $L_0\cap L_1=0$.

We now turn to families. Let $B$ be some manifold.
We consider a pair of smooth families of 
Lagrangian subspaces $B\ni b\mapsto L_i(b)$, $i=0,1$, and
we assume that $L_0(b)\cap L_1(b)=0$, $\forall b\in B$.
We obtain a corresponding family $\{\cD(b)\}_{b\in B}$ of invertible operators.
We want to apply the superconnection formalism in order to define
the $\eta$-form of that family. Since this formalism involves
derivatives of the family with respect to $b$ we prefer
to work with an unitarily equivalent family $\{\tilde{\cD}(b)\}_{b\in B}$ which has the advantage that its domain is independent of $b\in B$.

The $\eta$-form is a local object with respect to the base space.
In order to define it we only consider a germ of the family near a point $b_0\in B$.
Let $U_I(V)$ denote the group of unitary operators on $V$ 
which commute with $I$. The group $U_I(V)$ acts
transitively on the space $\Lambda$ of all Lagrangian subspaces of $V$.
Thus we can find germs of smooth 
families of unitaries $b\mapsto U_i(b)\in U_I(V)$, $i=0,1$,
with $U_i(b_0)=1$ and $U_i(b)L_i(b)=L_i(b_0)$.
We can define germs of smooth families $A_i(b):=\log(U_i(b))$ of anti-hermitean
matrices using the standard branch of the logarithm.
Let $\chi_i\in C^\infty([0,1])$ be cut-off functions with
$\chi_0(t)= 1$ for $t<1/5$, $\chi_0(t)= 0$ for $t>2/5$,
$\chi_1(t)=1$ for $t>4/5$, and $\chi_1(t)=0$ for $t<3/5$.
Then we set $W(t,b):=\exp(\chi_0(t) A_0(b)+\chi_1(t)A_1(b))$.
Then $b\mapsto W(.,b)$ can be considered as a germ of a family of unitary 
multiplication operators on $L^2([0,1],V)$.
We set 
\begin{eqnarray*}
\tilde{D}(b)&:=&W(.,b)DW^*(.,b)\\
&=& D- \chi_0^\prime I A_0(b) - \chi_1^\prime I A_1(b)\ ,
\end{eqnarray*}
where $"{}^\prime"$ denotes the derivative with respect to $t$.
We define the selfadjoint extension of $\tilde{\cD}(b)$ using the 
Lagrangian subspaces $L_0(b_0)$, $L_1(b_0)$.
Then $\tilde{\cD}(b)$ is unitarily equivalent to $\cD(b)$ and
its domain is independent of $b$.

We now turn to the definition of the $\eta$-form of the family $\{\tilde{\cD(b)}\}_{b\in B}$.
Let $C_1$ denote the graded algebra over $\C$ generated by $\sigma$
satisfying $\sigma^2=1$, $\sigma^*=\sigma$, and $\deg(\sigma)=1$. 
Let $\cH$ denote the germ at $b_0$ of the trivial Hilbert space bundle with fibre $L^2([0,1],V)\otimes C_1$ over $B$. 
We define the superconnection 
$A_s$, $s>0$, on $\cH$ associated to $\tilde{\cD}$ by 
$$A_s=\nabla-d(\chi_0 A_0+\chi_1 A_1)+\sqrt{s}\sigma\tilde{\cD}\ ,$$
where $d$ differentiates along $B$. 
Here $\nabla$ is the canonical connection of $\cH$ and $\nabla-d(\chi_0 A_0+\chi_1 A_1)=W\nabla W^*$.
For $\Ree(u)>1$ we can define the holomorphic family of germs of smooth, even
differential forms
\begin{equation}\label{forn1}\eta(u):=\frac{1}{2\sqrt{\pi}}\int_0^\infty \tr_\sigma^{even}(\sigma\tilde{ \cD}\ee^{-A_s^2}) s^{u-1/2}ds\ .\end{equation}
Here $\tr_\sigma^{even}(\dots)$ stands for the even form part of $\tr(\sigma\dots)$. 
As usual the asymptotic expansion of 
the heat kernels for small times implies
that $\eta(u)$ has a meromorphic continuation with respect to $u$
to all of $\C$ having at most first order poles.
Following \cite{daizhang96} we define
\begin{ddd}
$$\eta(\tilde{\cD}):= P.F.\: \eta(0)\ ,$$
\end{ddd}
where "$P.F.$" stands for the finite part of the Laurent expansion of $\eta(u)$
at $u=0$.
As defined above the form $\eta(\tilde{\cD})$ may depend on the choices
made for the definition of $\tilde{\cD}$. 
But the following lemma justifies the notation $\eta(L_0,L_1):=\eta(\tilde{\cD})$.
\begin{lem}
$\eta(\tilde{\cD})$ does not depend on the choices of 
the families $U_i$ and the cut-off functions $\chi_i$.
\end{lem}
\proof
Let $\hat{U}_i$, $\hat{\chi}_i$ be another choice and define
$\hat{W}$ as above. Let $\hat{\tilde{\cD}}$ denote the corresponding
family of operators. We set $V=W\hat{W}^*$. Then $\hat{\tilde{\cD}}=V^*\tilde{\cD}V$.
If $\hat{A}_s$ is the superconnection associated to $\hat{\tilde{\cD}}$,
then $\hat{A}_s=V^* A_s V$. It is now easy to check that
$\hat{\eta}(u)=\eta(u)$ if $\hat{\eta}(u)$ corresponds to $\hat{A}_s$.
\hB

\section{The Maslov cocycle}\label{sec1}

Now we turn to the generalized Maslov index.
Let $\{L_i(b)\}_{b\in B}$, $i=0,1,2$,
be smooth families of Lagrangian subspaces of $V$
such that $L_i(b)\cap L_j(b)=\{0\}$, $\forall b\in B$, $i\not= j$.
Let $\ch:K^0(B)\rightarrow H^{ev}(B,\R)$ be the Chern character.
In the present section we prove
\begin{prop}\label{sert}
The form $\eta(L_0,L_1)+\eta(L_1,L_2)+\eta(L_2,L_0)$
is closed. Moreover
$$\tau(L_0,L_1,L_2):=[\eta(L_0,L_1)+\eta(L_1,L_2)+\eta(L_2,L_0)]\in \ch(K^0(B))\subset H^{even}(B,\R)\ .$$
\end{prop}
\underline{Remark:} The zero component $\tau(L_0,L_1,L_2)^{0}\in \Z$
is the Maslov index (in its hermitean symplectic generalization) 
of the triple $(L_0,L_1,L_2)$ (see \cite{bunke95}). 
For an exposition of the usual Maslov index we refer to \cite{lionvergne80}).
In Proposition \ref{expl} below we explicitly compute the class $\tau(L_0,L_1,L_2)\in H^{even}(B,\R)$ in terms of the hermitean symplectic geometry of the family $\{L_0(b),L_1(b),L_2(b)\}_{b\in B}$.\\

\noindent
\proof
The idea of the proof is to formulate an index problem
for a family of Dirac operators $\dirac^+$ such that the form
$\eta(L_0,L_1)+\eta(L_1,L_2)+\eta(L_2,L_0)$
represents the Chern character of the index bundle of $\dirac^+$.

We consider a compact oriented Riemann surface $M_c$ which
is homeomorphic to a disc, and the boundary of which 
is the union of twelve pieces
$\partial_i M_c$, $i=0,\dots,11$. The boundary pieces
are labelled in their cyclic order.
We assume that $\partial_i M_c$ are isometric to the interval $[0,1]$, and that the boundary pieces
intersect in twelve corners. We assume that the metric is product in a neighbourhood of the interior of the
pieces and that neighbourhoods of 
the corners are isometric to a neighbourhood of the vertex of the
euclidean quadrant $\R^+\times\R^+$.

Let $M$ be the oriented, non-compact Riemann surface with six boundary components $\partial_i M$, $i=0,\dots, 5$,  isomorphic to $\R$ 
which is obtained by gluing infinite cylinders $[0,\infty)\times [0,1]$ along
the boundary pieces of $M_c$ with odd label.
The boundary components of $M$ are again labelled according to their cyclic order.

We consider the spinor bundle $S=S^+\oplus S^-$ of $M$ 
and fix a finite-dimensional Hilbert space $W$ of dimension $\dim(V)/2$.
We consider the graded vector bundle $E=S\otimes (W\oplus W^{op})$.
By $\dirac$ we denote the corresponding twisted Dirac operator.
Let $\dirac^\pm$ be the parts
mapping sections of $E^\pm$ to those of $E^\mp$.
We formulate an index problem for $\dirac^+$ by putting
boundary conditions at $\partial_j M$ depending on $L_i(b)$, $i=0,1,2$.

The metric of $M$ is flat near infinity and the boundaries (the flat region).
We claim that the holonomy of the parallel transport in $S$ along $\partial M_c$
is trivial, where we consider $M_c$ as a submanifold of $M$. 
Note that $M$ is topologically a disc.
Thus we can choose a trivialization of the tangent bundle $TM$.
Measured with respect to the trivialization the parallel transport in $TM$ along $\partial M_c$ gives a rotation of $-4\pi$. 
The trivialization of $TM$ induces one of the $Spin(2)$-principal
bundle of $M$. The parallel transport along $\partial M_c$
in this bundle corresponds to a rotation of $-2\pi$ in the
structure group. This implies the claim.
We fix some point in $\partial M_c$ and identify the fibres of the
bundle $S$ near infinity and $\partial M$ 
with the fibre over this point using the parallel transport.
We denote this fibre by $S$, too.
Analogously we also denote the fibre of $E$ over this point
by $E$. 

Let $(s,t)$ be oriented euclidean orthonormal 
coordinates near a point in the flat region.
Then we have 
$\dirac=\sigma_s \partial_s + \sigma_t\partial_t$, where
$\sigma_s,\sigma_t\in \Hom^{odd}(E,E)$ depend on the choice of coordinates.
Again we consider the components $\sigma^\pm_s,\sigma^\pm_t\in \Hom^{odd}(E^\pm,E^\mp)$.
The operator $I:=\sigma_s^-\sigma_t^+\in\Aut(E^+)$ is invariantly defined.
It satisfies $I^*=-1$, $I^2=-1$, $\tr I=0$, and it defines a hermitean
symplectic structure on $E^+$. We fix an isometry $V\cong E^+$ which is compatible with the complex structures $I$ on $V$ and $E^+$.

Now we introduce the family of boundary conditions defining the family
$\{\dirac^+(b)\}_{b\in B}$.
Let $C^\infty_{L^2}(M,E^\pm)$ denote the space of all
sections of $E^\pm$ which are square integrable together with all their derivatives.
We let $\dirac^+(b)$ be the differential operator $\dirac^+$ mapping 
$$\{\phi\in C_{L^2}^\infty(M,E^+)\:|\: \phi(x)\in L_i(b)\forall x\in\partial_i M\:\mbox{or}\: x\in\partial_{i+3} M\}$$ to
$C_{L^2}^\infty(M,E^-)$.

First we show that $\dirac^+(b)$ gives rise to a smooth family
of Fredholm operators such that the index bundle is well-defined.
Then we apply the superconnection formalism in order to compute the Chern character of the index bundle.

First we conjugate the family $\{\dirac^+(b)\}_{b\in B}$ to a family $\{\tilde{\dirac}^+(b)\}_{b\in B}$ with constant domain.
This will be done again on the level of germs at a point $b_0\in B$.
In the remainder of the present section we replace $B$ by a sufficiently small neighbourhood of $b_0$.
We define a germ of a family of smooth $U(V)$-valued functions
$W^+(b,m)$, $m\in M$, such that near $\partial_i M$, $\partial_{i+3} M$,
\begin{equation}\label{eqq1}W^+(b,(s,t))=\exp(\chi_0(s)A_i(b))\ ,\end{equation}
where $A_i$, $\chi_0$ were defined above.
Here $(s,t)$ are orthonormal euclidean coordinates, $s$ being normal
to the boundary given by $s=0$. 
$W^+$ is determined by (\ref{eqq1}) near $\partial M$ and we continue $W^+$
to the interior of $M$ by the constant $1\in U(V)$. Similarly we set $W^-=-\sigma_s^+ W^+\sigma_s^-$ near $\partial M$ and
continue $W^-$ to the interior of $M$ by $1\in E^-$.

Let $\tilde{\dirac}^+:= W^-\dirac (W^+)^*$. Then $\tilde{\dirac}^+(b)$
is a germ of a family of (now $b$-dependent) Dirac operators with domain
(now independent of $b$) given as above by the Lagrangian subspaces
$L_i(b_0)$ at $\partial_i M$, $\partial_{i+3} M$.

For $i=0,\dots,5$ the family $\tilde{\dirac}^+(b)$ induces families of translation
invariant operators $\tilde{\dirac}_i^+(b)$ on the infinite cylinder
$\R\times [0,1] $. Let $\tilde{\dirac}_i^-(b)$ be the
formal adjoint (with the adjoint boundary condition) of $\tilde{\dirac}_i^+(b)$
and let 
$$\tilde{\dirac}_i(b):=\left(\begin{array}{cc}0&\tilde{\dirac}_i^-(b)\\ \tilde{\dirac}_i^+(b)&0\end{array}\right)\ .$$

It is easy to check that the boundary condition
satisfies the condition of  Lopatinski-Shapiro.
Hence $\tilde{\dirac}_i(b)$ has a parametrix of bounded propagation.
It follows that $\tilde{\dirac}_i(b)$ is essentially selfadjoint
on the space of smooth sections with compact support satisfying
the boundary condition. Using separation of variables and
the transversality of the Lagrangian subspaces one can then
check that $\tilde{\dirac}_i(b)$ is invertible.

Gluing the distributional kernels of the $(\tilde{\dirac}_i(b))^{-1}$, $i=0,\dots,5$, 
at infinity with a local parametrix of $\tilde{\dirac}(b)$
we obtain a parametrix $Q(b)$ of $\tilde{\dirac}(b)$ depending smoothly
on $b$. In particular, the smoothing remainders $\tilde{\dirac}(b) Q(b)-1$, $Q(b)\tilde{\dirac}(b) -1$ have compact support. 
$\{\tilde{\dirac}(b)\}$ is a family of essentially selfadjoint
operators. The domain of definition $\cH$ of the unique selfadjoint extension
of $\tilde{\dirac}(b)$ is independent of $b$. 
When viewed as a family of bounded operators
from $\cH^+$ to $L^2(M,E^-)$ the family $\{\tilde{\dirac}^+(b)\}$
is a smooth family of Fredholm operators and the index bundle of $\tilde{\dirac}^+(b)$ is well defined.

We now apply the superconnection formalism 
in order to obtain a formula for the Chern
character $\ch(\ind(\tilde{\dirac}^+))$ of the index
bundle of $\{\tilde{\dirac}^+(b)\}_{b\in B}$. 
Our main goal is to work out the computation which
gives the $\eta$-form as the boundary contribution to 
Chern character of the index bundle. 
We supress the standard arguments (see \cite{berlinegetzlervergne92}, Ch.9.5, \cite{melrosepiazza96}) which show that
the small time analysis of the heat trace indeed gives a representative of the
Chern form of the index bundle.
In particular, we assume for simplicity that $\dim\ker\tilde{\dirac}^+(b)$
is constant.

Let $\nabla $ denote the trivial connection on the 
bundle $B\times \cH$. 
We set $W:=W^+\oplus W^-$.
Let $\tilde{\nabla}=W\nabla W^*$.
Then we define the superconnection
$$B_t=\tilde{\nabla}+\sqrt{t}\tilde{\dirac}\ .$$
The curvature $B_t^2$ has the form
$$B_t^2=t\tilde{\dirac}^2+\sqrt{t} R\ ,$$
where $R$ is a one-form with values in the odd
endomorphisms of $E$.
The heat operator $\ee^{-B_t^2}$ can be constructed
using the Volterra series \cite{berlinegetzlervergne92}, Prop.9.46.
Let $P_t(x,y)$ denote the smooth integral kernel of $\ee^{-B_t^2}$.
In order to simplify the notation we omit the smooth dependence on $b\in B$.

Note that $M$ is a manifold with a cylindrical end
$N\times [0,\infty)$, where $N$ is isometric to the disjoint union
of six copies of the unit interval.
Let $(n,r)$ denote corresponding coordinates. 
Though $\ee^{-B_t^2}$ is not of trace class we define
$$\Tr^\prime_s \ee^{-B_t^2}:=
\int_{M_c} \tr_s P_t(x,x) dx +
\lim_{u\to\infty} \int_0^u \int_{N} \tr_s P_t((n,r),(n,r)) dn \:dr\ .$$

We first argue that this limit exists.
We claim that for some $C<\infty$, $c>0$,
$$|\tr_s P_t((n,r),(n,r))|<C \ee^{-cr^2}\ .$$
The constants $C,c$ can be choosen uniformly for
$t$ varying in compact subsets of $(0,\infty)$.
Consider the infinite cylinder $Z:=N\times \R$.
Let $E^Z$ be the bundle on $Z$ induced by $E$.
Let $\{\tilde{\dirac}^Z(b)\}_{b\in B}$ denote the
family of translation invariant operators on $E^Z$ induced by $\{\tilde{\dirac}(b)\}_{b\in B}$. Let $\cH^Z$ the domain
of $\tilde{\dirac}(b)^Z$ which is again independent of $b\in B$. 
We then obtain a translation invariant  
superconnection $B^Z_t$ on the bundle $B\times \cH^Z$. 
Let $P^Z_t((n,r),(m,s))$ denote the corresponding heat kernel.
The Clifford multiplication by $\imath\sigma_r$ is unitary, odd, and commutes with  $(B_t^Z)^2$. Thus
$$\tr_s P^Z_t((n,r),(n,r))\equiv 0\ .$$
Using standard finite-propagation speed estimates one can show that for $r>0$
$$|P^Z_t((n,r),(n,r))-P_t((n,r),(n,r))|<C\ee^{-cr^2}\ ,$$
and this proves the claim.

Let $\tilde{\cD}^Z$ be the family of operators on $N\times V$ induced by
$\tilde{\dirac}^Z$. Then $\tilde{\cD}^Z$ can be identified with the dirct sum of two copies of the direct sum of three copies of $\tilde{\cD}$ with boundary conditions given by the families of pairs $(L_0,L_1)$, $(L_1,L_2)$, $(L_2,L_0)$.
By $A_t^Z$ we denote the superconnection corresponding to $\tilde{\cD}^Z$. 

Set $\gamma:=\sigma_r\in \End(E)$.
The comparison with the cylinder $Z$
shows that $\Tr^\prime_s \ee^{-B_t^2}$ can be differentiated with respect to $t$. Using Duhamel's formula one can check that one can
commute $\Tr_s^\prime$ and $d/dt$.
 \begin{lem}
$$\frac{d}{dt} \Tr^\prime_s \ee^{-B_t^2}= -\frac{1}{\sqrt{4\pi t}} \Tr_s
\gamma \tilde{\cD} \ee^{-(A_t^Z)^2}-\frac{1}{2\sqrt{t}} d\Tr^\prime_s  \tilde{\dirac}\ee^{-B^2_t}$$
\end{lem}
\proof
First we claim that
$$\frac{d}{dt} \Tr^\prime_s \ee^{-B_t^2}=-\Tr^\prime_s [B_t,\frac{dB_t}{dt}\ee^{B^2_t}]\ .$$
Let $\rho_u$ denote the characteristic function of $M_c\cup N\times [0,u]$.
Using Duhamel's formula we get
\begin{eqnarray*}
\frac{d}{dt} \Tr^\prime_s \ee^{-B_t^2}&=&
-\Tr^\prime_s \int_0^1 \ee^{-s B^2_t}\frac{dB^2_t}{dt}\ee^{-(1-s)B^2_t} ds\\
&=&-\lim_{u\to\infty} \int_0^1  \Tr_s \rho_u \ee^{-s B^2_t}\frac{dB^2_t}{dt}\ee^{-(1-s)B^2_t} ds\\
&=&-\lim_{u\to\infty} \int_0^1 \Tr_s[\frac{dB_t}{dt},B_t] \ee^{-s B^2_t}\rho_u\ee^{-(1-s)B^2_t} ds\\
&=&-\lim_{u\to\infty} \lim_{v\to\infty} \int_0^1 \Tr_s \rho_v [\frac{dB_t}{dt},B_t] \ee^{-s B^2_t}\rho_u\ee^{-(1-s)B^2_t}ds\\
&=&-\lim_{v\to\infty} \lim_{u\to\infty}  \int_0^1 \Tr_s \rho_v [\frac{dB_t}{dt},B_t] \ee^{-s B^2_t}\rho_u\ee^{-(1-s)B^2_t}ds\\
&=&-\lim_{v\to\infty}  \int_0^1 \Tr_s \rho_v [\frac{dB_t}{dt},B_t] \ee^{- B^2_t} ds \\
&=&-\Tr^\prime_s [B_t,\frac{dB_t}{dt}\ee^{-B^2_t}]\ .\end{eqnarray*}
In order to justify that $\lim_{v\to\infty}$ and $\lim_{u\to\infty}$
can be interchanged
one can again use the comparison with the infinite cylinder $Z$.
We use $$\frac{dB_t}{dt}=\frac{1}{2\sqrt{t}}\tilde{\dirac}$$ in order to write
\begin{eqnarray}
-\Tr^\prime_s [B_t,\frac{dB_t}{dt}\ee^{B^2_t}]&=&
-\frac{1}{2\sqrt{t}}\Tr^\prime_s [\tilde{\nabla},\tilde{\dirac}\ee^{-B^2_t}] -
\frac{1}{2}\Tr_s^\prime [\tilde{\dirac},\tilde{\dirac}\ee^{-B^2_t}]\label{qq1}\\
&=&-\frac{1}{2\sqrt{t}} d\Tr^\prime_s  \tilde{\dirac}\ee^{-B^2_t} -\frac{1}{2}\Tr_s^\prime [\tilde{\dirac},\tilde{\dirac}\ee^{-B^2_t}]\label{q1}\ .
\end{eqnarray}
Before justifying the transition from (\ref{qq1}) to (\ref{q1}) we
consider the second term of (\ref{q1}).
Let $z$ denote the $\Z_2$ grading operator.
By integration by parts we obtain
\begin{eqnarray*}
\lefteqn{-\frac{1}{2}\Tr_s^\prime [\tilde{\dirac},\tilde{\dirac}\ee^{-B^2_t}]}\hspace{0.5cm}\\
&=&-\frac{1}{2}\Tr^\prime z \tilde{\dirac}^2\ee^{-B^2_t}-\frac{1}{2}\Tr^\prime z \tilde{\dirac}\ee^{-B^2_t}\tilde{\dirac}\\
&=&-\frac{1}{2}\lim_{u\to\infty} \int_0^u \int_{N}\tr z (\tilde{\dirac}^2 \ee^{-B_t})((n,u),(n,u))dndu-\frac{1}{2}\int_{M_c}\tr z (\tilde{\dirac}^2\ee^{-B^2_t})(m,m) dm\\
&&-\frac{1}{2}\lim_{u\to\infty} \int_0^u \int_{N}\tr z (\tilde{\dirac} \ee^{-B_t}\tilde{\dirac})((n,u),(n,u))dndu-\frac{1}{2}\int_{M_c}\tr z (\tilde{\dirac}\ee^{-B^2_t}\tilde{\dirac})(m,m) dm\\
&=&-\frac{1}{2}\lim_{u\to\infty} \int_{N} \tr (\gamma z\tilde{\dirac} \ee^{-B_t} )((n,u),(n,u))\\
&=&\frac{1}{2}\lim_{u\to\infty} \int_{N} \tr_s (\gamma \tilde{\dirac} \ee^{-B_t} )((n,u),(n,u))\ .
\end{eqnarray*}
In order to evaluate this limit we can replace the kernel $P_t$ by $P_t^Z$.
We use the Volterra series \cite{berlinegetzlervergne92}, Prop.9.46, in order to compute $P_t^Z$.
Note that on $Z$
\begin{equation}\label{q2}\ee^{-t(\tilde{\dirac}^Z)^2}((n,r),(m,s))=\ee^{-t(\tilde{\cD}^Z)^2}(n,m)\frac{\ee^{-(r-s)^2/4t}}{\sqrt{4\pi t}}\ ,\end{equation}
and its follows
$$P^Z_t=\ee^{-t(\tilde{\dirac}^Z)^2}+\sum_{k=1}^\infty (-1)^k t^{k/2} 
\int_{\Delta^k}\ee^{-t\sigma_0(\tilde{\dirac}^Z)^2} R^Z\dots R^Z\ee^{-t\sigma_k(\tilde{\dirac}^Z)^2} d\sigma\ ,$$
where $\Delta_k$ denotes the standard $k$-simplex.
Inserting (\ref{q2}) we obtain
\begin{eqnarray*} \lefteqn{P^Z_t((n,r),(m,s))}\hspace{0.3cm}\\
&=&\frac{\ee^{-(r-s)^2/4t}}{\sqrt{4\pi t}}\left( \ee^{-t(\tilde{\cD}^Z)^2}+\sum_{k=1}^\infty (-1)^k t^{k/2} 
\int_{\Delta^k}\ee^{-t\sigma_0(\tilde{\cD}^Z)^2} R^Z\dots R^Z\ee^{-t\sigma_k(\tilde{\cD}^Z)^2} d\sigma \right)(n,m)\ .\end{eqnarray*}
We now apply $$\gamma\tilde{\dirac}^Z=-\frac{\partial}{\partial r}+\gamma\tilde{\cD}^Z$$
and evaluate the result at $r=s$.
We then obtain
$$\frac{1}{2}\lim_{u\to\infty} \int_{N} \tr_s\gamma (\tilde{\dirac} \ee^{-B^2_t})(n,u) dn = \frac{1}{2}\lim_{u\to\infty} \int_{N} \tr_s\gamma (\tilde{\dirac}^Z \ee^{-(B^Z_t)^2})(n,u) dn = \frac{1}{4\sqrt{\pi t}} \Tr_s \gamma \tilde{\cD}^Z\ee^{-(A^Z_t)^2}\ . $$
Now we consider the first term of (\ref{q1}).
By a similar computation as above one can show that
on the cylinder $Z$
$$\tr_s [\tilde{\nabla}^Z,\tilde{\dirac}^Z \ee^{-(B^Z_t)^2}]((n,r),(n,r))\equiv 0\ .$$
Thus on $M$ this quantity vanishes rapidly as $r\to\infty$. 
We can take $\Tr_s^\prime$ and
$$-\frac{1}{2\sqrt{t}}\Tr^\prime_s [\tilde{\nabla},\tilde{\dirac} \ee^{-B^2_t}]= -\frac{1}{2\sqrt{t}} d\Tr^\prime_s  \tilde{\dirac}\ee^{-B^2_t}$$
is an exact form.
This finishes the proof of the lemma. \hB

Let $\nabla^0$ denote the induced connection on the index bundle of $\dirac^+$
and let $\ch(\nabla^0)=\tr_s\ee^{-(\nabla^0)^2}$ be the corresponding Chern form.
As in \cite{berlinegetzlervergne92}, Ch. 9 or \cite{melrosepiazza96}
one can show the following estimates.
\begin{lem}
Let $|.|$ be any continuous seminorm on the space of smooth
forms on $B$.
For $t\to\infty$ we have
\begin{eqnarray*}
|\Tr^\prime_s \ee^{-B_t^2}-\ch(\nabla^0)|&=&O(t^{-1/2})\\
|\frac{1}{2\sqrt{t}} \Tr^\prime_s  \tilde{\dirac}\ee^{-B^2_t}|&=& O(t^{-3/2})\ .
\end{eqnarray*}
\end{lem}
We conclude that
$$ \int_s^\infty \frac{1}{4\sqrt{\pi t}} \Tr_s
\gamma \tilde{\cD^Z} \ee^{-(A_t^Z)^2} dt=:\hat{\eta}(s)$$ exists and that
\begin{equation}\label{sab}\ch(\nabla^0)=\Tr^\prime_s \ee^{-B_s^2}+\hat{\eta}(s) + d \alpha(s)\ ,\end{equation}
where
\begin{equation}\label{r1}\alpha(s):=-\int_s^\infty \frac{1}{2\sqrt{t}} \Tr^\prime_s  \tilde{\dirac}\ee^{-B^2_t} dt\ .\end{equation}
We now want to take the limit $s\to 0$.
\begin{lem}
Let $|.|$ be any continuous seminorm on the space of smooth
forms on $B$.
For $t\to 0$ we have
\begin{eqnarray}
|\Tr_s^\prime \ee^{-B^2}|&=&o(1)\label{w1}\\
|\frac{1}{2\sqrt{t}} \Tr^\prime_s  \tilde{\dirac}\ee^{-B^2_t}|&=& O(t^{-1/2})\ .\label{w2}
\end{eqnarray}
\end{lem}
\proof 
We first show (\ref{w1}). We employ finite propagation speed estimates and the comparison with the cylinder
in order to show that on the end of $M$
$$|\tr_s P_t (n,r)|<C\ee^{-r^2/t}\ .$$
The local index theorem \cite{berlinegetzlervergne92}, Ch.10,
gives $|\tr_s P_t (x)|=o(1)$ for $x\in M$ since $S$ is twisted with a bundle
of the form $W\oplus W^{op}$. Both estimates together give (\ref{w1}).

Now we consider (\ref{w2}). On the cylinder $Z$
we have $\tr_s \tilde{\dirac}^Z \ee^{-(B_t^Z)^2}(n,r)\equiv 0$.
We conclude that on $N\times[0,\infty)\subset M$ 
$|\tr_s \tilde{\dirac} \ee^{-B^2_t} (n,r)|<C\ee^{-r^2/t}$.
Moreover, for $x$ in a small neighbourhood of the boundary of $M_c$ we
have $|\tr_s \tilde{\dirac} \ee^{-B^2_t} (x)|<C\ee^{-c/t}$.
If $x$ is in the interior of $M_c$ we can employ the method of \cite{berlinegetzlervergne92}, Ch. 10.5, in order to show that 
$$|\frac{1}{2\sqrt{t}} \tr_s \tilde{\dirac} \ee^{-B^2_t} (x)|=O(t^{-1/2})\ .$$
The estimate (\ref{w2}) follows.
\hB
Now we can take the limit $s\to 0$ in (\ref{r1}).
We obtain 
$$\ch(\nabla^0)=\hat{\eta}(0)+d\alpha(0)\ .$$
If $\eta(u,L_i,L_j)$ denotes the form (\ref{forn1})
which is defined using the boundary condition given by family of pairs $(L_i,L_j)$,
then $\hat{\eta}(0)=\lim_{u\to 0}(\eta(u,L_0,L_1)+\eta(u,L_1,L_2)+\eta(u,L_2,L_0))$.
We conclude that
$$\res_{u=0}(\eta(u,L_0,L_1)+\eta(u,L_1,L_2)+\eta(u,L_2,L_0))=0$$
and  that $\tau(L_0,L_1,L_2)=\ch(\nabla^0)+d\alpha(0)$.
This proves Proposition \ref{sert}.\hB

\section{Computation of $\tau(L_0,L_1,L_2)$}
 
Let $L_0,L_1,L_2\subset V$ be pairwise transverse Lagrangian subspaces.
 Then $V=L_1\oplus L_2$ and we can write $x_0=x_1+x_2$,
$x_i\in L_i$, $i=0,1,2$.  We define a hermitean quadratic form $Q$ on $L_0$ by
$$Q(x_0):=(Ix_1,x_2)\ ,$$
where $(I,.,)$ is the symplectic form of $V$ associated to $I$ and the Hilbert
space structure of $V$. It is easy to see that $Q$ is nondegenerate.
Thus we can split $L_0=L_0^+\oplus L_0^-$ into the positive and negative
eigenspace of $Q$.
Returning now to the family case we obtain a decomposition 
$L_0=L_0^+\oplus L_0^-$ of the bundle
of Lagrangian subspaces $L_0\subset B\times V$ which is induced by
the two other subbundles $L_1,L_2$.
\begin{prop}\label{expl}
We have
$$\tau(L_0,L_1,L_2)=\ch(L_0^+)-\ch(L_0^-)\in H^*(B,\Z)\ .$$
\end{prop}

\proof
The proof of the proposition consists of two steps.
\begin{enumerate}
\item Using the $K$-theoretic relative index theorem
\cite{bunke9381} we reduce to an index problem for a family
of Dirac operators on the disc. The parameter dependence of this family is again built in through the boundary conditions.
\item We then consider the "universal" family of such operators
which is parametrized by a space which is homotopy equivalent to the space of all triples of pairwise transverse Lagrangian subspaces of $V$. It suffices to verify
the assertion of the proposition in this special case.
\end{enumerate}
First we want to compactify $M$ by cutting off the cylindrical ends
and gluing in half discs. The resulting manifold $\hat{M}$ is
then topologically a disc. 
Let $\hat{\dirac}^+$ be the corresponding Dirac operator. We want to find
a family of boundary conditions parametrized by $B$ such that
$\ch(\ind(\dirac^+))=\ch(\ind(\hat{\dirac}^+))$.

Let $Y\subset \R^2$ denote the subset
$$Y:=\{(s,t)\in\R^2\:|\:s\ge 0, t\in [-1/2,1/2]\:\mbox{or}\: s\le 0, t^2+s^2\le 1/4\}\ .$$
Then $Y$ is a Riemannian surface with $C^1$-boundary and one cylindrical end.
Let $S_Y$ be the spinor bundle of $Y$.
Let $E_Y:=S_Y\otimes (W\oplus W^{op})$ and $\dirac_Y$ be the Dirac operator on $E_Y$.
We trivialize $S_Y$ and $E_Y$ using the flat Levi Civita connection and denote the
typical fibre of $E_Y$ by $E$.
Then $$\dirac_Y^\pm=\sigma_s^\pm \frac{\partial}{\partial s} + \sigma^\pm_t \frac{\partial}{\partial t}\ .$$
The space $V:=E^+$ is a symplectic vector space with symplectic structure induced by $I=\sigma_s^-\sigma_t^+$.

Let now $L_0,L_1$ be transversal Lagrangian subspaces of $V$. We want to construct a 
family of Lagrangian subspaces $L(s,t)=L(s,t)(L_0,L_1)$ which is parametrized by $(s,t)\in\partial Y$, and which depends smoothly on $L_0,L_1$, such that
$L(s,t)=L_0$ for $s\ge 0$, $t=1/2$ and $L(s,t)=L_1$ for $s\ge 1$, $t=-1/2$.
If $B\ni b\rightarrow (L_0(b),L_1(b))$ is a smooth family of pairs of transverse Lagrangian subspaces,
then we require that the index of the Dirac operator $\dirac^+_Y$ subject to the family of boundary
conditions $L(.,.)(L_0(b),L_1(b))$ is trivial.

We first set 
$$\hat{L}(s,t)(L_0):=\left\{\begin{array}{cc} L_0 & (s,t)\in\partial Y, t\ge 0\\
                                           IL_0 & (s,t)\in \partial Y, s\ge 0, t=-1/2\\
                                            \sigma^-_sn(s,t) L_0&  (s,t)\in \partial Y, s\le 0, t\le 0
\end{array}\right. \ ,$$
where $n(s,t):=2(s\sigma_s^++t\sigma_t^+)$ is the Clifford multiplication by the normal vector. 
It is easy to see that $\hat{L}(s,t)$ is $C^1$ with respect to $(s,t)$.
If $B\ni b\rightarrow L_0(b)$ is a smooth family of Lagrangian subspaces,
then we consider the family $\{\hat{\dirac}_Y^+(b)\}_{b\in B}$ given by $\dirac_Y^+$ subject to the boundary conditions
given by $\hat{L}(.,.)(L_0(b))$.     
\begin{lem}
 In $K^0(B)$ we have $2\: \ind(\hat{\dirac}_Y^+)=0$.
\end{lem}
\proof We claim that $\hat{\dirac}_Y^+$ is equivalent with its adjoint $(\hat{\dirac}_Y^+)^*$.
Thus $2\:\ind (\hat{\dirac}_Y^+)=\ind(\hat{\dirac}_Y^+)-\ind(\hat{\dirac}_Y^+)^*=0$.

We now show the claim. First we describe the adjoint $(\hat{\dirac}^+)^*$.
Note that $E^-$ is a hermitean symplectic vector space 
with symplectic structure induced by $\sigma^+_s\sigma^-_t$.
Then  the adjoint of $\hat{\dirac}_Y^+$ is the operator $\hat{\dirac}_Y^-$ subject to the boundary conditions given
by the family of Lagrangian subspaces of $E^-$
$$B\ni b\rightarrow \{(s,t)\in\partial Y\mapsto \tau(s,t)\hat{L}(s,t)(L_0(b))\}\ ,$$ where $\tau(s,t)$
is the Clifford multiplication by the tangent vector at $(s,t)\in \partial Y$.
If $\psi^+(s,t)$ is a section of $E^+_Y$, then we set $(U^+\psi^+)(s,t):=\sigma_t^+\psi^+(s,-t)$. Then $U^+$ is unitary,
and $U^+\psi^+$ is a section of $E_Y^-$. It is easy to check that
the equality of differential operators
$$U^+\hat{\dirac}_Y^-U^+=\hat{\dirac}_Y^+$$
is compatible with the boundary conditions. This shows the claim. \hB

Let $\Omega(x,y)=(Ix,y)$ be the (hermitean) symplectic form on $V$.
Let $\Lambda$ denote the manifold of all Lagrangian subspaces of $V$.
For $L\in\Lambda$ let $\cL_L$ denote the subset of all Lagrangian subspaces $L^\prime\in\Lambda$
which are transverse to $L$. The following discussion is parallel to that in \cite{guilleminsternberg77} p. 117/118. Let $P_{L^\prime}$ denote the projection from $V$ to $L^\prime$ with kernel $L$.
It is easy to check that $\Omega(P_{L^\prime}x,y)+\Omega(x,P_{L^\prime}y)=\Omega(x,y)$.
We define the hermitean quadratic form 
\begin{equation}\label{yu2}Q_{L^\prime}(x,y):=\Omega(P_{L^\prime}x,y)-\frac{1}{2}\Omega(x,y)\ .\end{equation}
Indeed
\begin{eqnarray*}
Q_{L^\prime}(x,y)&=&\Omega(P_{L^\prime}x,y)-\frac{1}{2}\Omega(x,y)\\
&=&-\Omega(x,P_{L^\prime}y)+\frac{1}{2}\Omega(x,y)\\
&=&\overline{\Omega(P_{L^\prime}y,x)-\frac{1}{2}\Omega(y,x)}\\
&=&\overline{Q_{L^\prime}(y,x)}\ .
\end{eqnarray*}
We have 
\begin{equation}\label{yu1}Q_{L^\prime}(x,y)=-\frac{1}{2}\Omega(x,y),\quad\forall x\in L, y\in V\ .\end{equation}
Any  hermitean quadratic form $Q$ satisfying (\ref{yu1}) determines a Lagrangian subspace
$L^\prime$ such that $Q=Q_{L^\prime}$. In fact let $P^\prime$ be determined by $Q$ and (\ref{yu2}), then
$L^\prime$ is just the $1$-eigenspace of $P^\prime$.
Thus we can identify $\cL_L$ with the space of hermitean quadratic forms satisfying (\ref{yu1}).
In particular, $\cL_L$ is an affine space where the affine structure only depends on $L$.

We now can construct the desired family $L(s,t)(L_0,L_1)$. 
Note that $IL_0,L_1\in \cL_{L_0}$, and there is a natural affine path
$L(r)=L(r)(L_0,L_1)$ with $L(0)=IL_0$, $L(1)=L_1$. We choose a smooth cut-off function
$\chi\in C^\infty([0,1])$ with $\chi(t)\in[0,1]$, $\chi(t)=0$ near $t=0$ and
$\chi(t)=1$ near $t=1$.
We set $L(s,t)=\hat{L}(s,t)$ for all $(s,t)\in\partial Y$ except
for $t=-1/2$, $s\in[0,1]$, where we set $L(s,t)=L(\chi(s))(L_0,L_1)$.
Then $L(s,t)$ is $C^1$ with respect to $(s,t)$ and depends smoothly on $L_0,L_1$.  Let $\{\dirac_Y^+(b)\}_{b\in B}$ denote the family of Dirac operators
given by $\dirac_Y^+$ subject to the boundary conditions defined by $\{L(.,.)(L_0(b),L_1(b))\}$. Then $\hat{\dirac}^+_Y$ and $\dirac_Y^+$ are
homotopic families and thus $2\:\ind(\dirac_Y^+)=0$ in $K^0(B)$.
The upshot of the construction above is that we associated to a pair
$L_0,L_1$ of transversal Lagrangian subspaces a canonical path $\gamma(L_0,L_1)$
from $L_0$ to $L_1$ which is parametrized by $\partial Y$.
The path $\gamma(L_0,L_1)$ depends smoothly on the pair $(L_0,L_1)$
and has (in a certain sense that will become clear below) the minimal winding number.

Now we can cut-off the six cylindrical ends of $M$ and glue in the pieces
$Y_i=\{(s,t)\in Y\:|\: s\ge 1\}$, $i=0,\dots,5$.
The resulting manifold $\hat{M}$ is topologically a two disc.

To be more precise let $Z_i=[0,\infty)\times [-1/2,1/2]$, $i=0,\dots,5$,
denote the cylindrical ends of $M$. Then we cut at
$\{1\}\times [-1/2,1/2]$. We identify
$(s,t)\in Z_i$ with $(2-s,t)\in Y_i$. Moreover,
we use $\sigma_s:E_{|\partial Z_i}\rightarrow E_{Y_i|\partial Y_i}$ in order to glue the bundles.
Then $\dirac^+$ glues with $\dirac_Y^-$.
Assume that on the component $[0,\infty)\times \{-1/2\}$ of $\partial Z_i$
we have the boundary condition given by $L_i$ (resp. $L_{i-3}$) and on 
$[0,\infty)\times \{1/2\}$ we have the one given by $L_{i+1}$ (resp. $L_{i-2}$), where $L_3=L_0$.
Then on the boundary part of $\hat{M}$ which comes from $Y_i$ we choose 
the path $\sigma_s\tau\gamma(L_i,L_{i+1})$ (resp. $\sigma_s\tau\gamma(L_{i-3},L_{i-2})$), where $\tau$ again denotes the Clifford
multiplication with the unit vector tangent to the boundary.
This path indeed connects $L_i$ with $L_{i+1}$. 
Thus we have constructed a closed path $\hat{\gamma}(L_0,L_1,L_2)$ of Lagrangian
subspaces of $V$ which is parametrized by $\partial\hat{M}$, and which depends smoothly on the triple $(L_0,L_1,L_2)$.

We use this path in order to define the boundary condition   
for the $W\oplus W^{op}$-twisted Dirac operator $\hat{\dirac}^+$ on $\hat{M}$.  
Recall that we identify $V$ with the fibres of the bundle $E^+_{\partial \hat{M}}$ using the parallel transport along $\partial \hat{M}$.
Using the $K$-theoretic relative index theorem \cite{bunke9381} it is easy to see that $\ind(\dirac^+)=\ind(\hat{\dirac}^+)\in K^0(B)[1/2]$.
Indeed, the relative index theorem states that
$$\ind(\dirac^+)+\sum_{i=0}^5\ind (\dirac_{Y_i}^+) = \ind(\hat{\dirac}^+)+ \sum_{i=0}^5\ind(\dirac_{Z_i}^+)\ ,$$
where $Z_i=[-1/2,1/2]\times \R$ and $\dirac_{Z_i}^+$ is the $W\oplus W^{op}$-twisted Dirac operator subject to the boundary conditions
given by $L_i$ at $\{-1/2\}\times \R$, $L_{i+1}$ at $\{-1/2\}\times \R$ 
(resp. $L_{i-3}$ at $\{-1/2\}\times \R$, $L_{i-2}$ at $\{-1/2\}\times \R$).
But $\ind(\dirac_{Z_i}^+)=0$ in $K^0(B)[1/2]$ for symmetry reasons.

Deforming the metric of $\hat{M}$ to the standard metric of the two disc
we do not change the index. Below we will assume that $\hat{M}$ is isometric
to the two disc. The parallel transport in $E^+$ along $\partial \hat{M}$
with respect to the globally flat metric gives an 
identification of $V$ with the fibres of $E^+$ which is topologically different from the one used above. This fact has to be taken into account below.
We have now finished the first part of the proof of the proposition.

We start with the second part. Let $\Lambda^3$ be the space of triples
$(L_0,L_1,L_2)$ of pairwise transverse Lagrangian subspaces of $V$.
Let $Sp(V)$ denote the group of symplectic automorphisms of $V$.
Note that $i\Omega$ is a non-degenerate hermitean form of
signature $(l,l)$, where $l=\dim_\C(V)/2$. Thus $Sp(V)\cong U(l,l)$.
The group $Sp(V)$ acts on $\Lambda^3$. We claim that $\Lambda^3$
is the disjoint union of orbits of $Sp(V)$.

First it is easy to see that $Sp(V)$ acts transitively on the space $\Lambda$.
Let $L_0\in \Lambda$. Then any $L_1\in \cL_{L_0}$ can be written
as $\{Bx+x\:|\: x\in IL_0\}$ for some $B\in \End(IL_0,L_0)$.
The condition that $L_1$ is Lagrangian translates to
$\Omega(Bx,y)+\Omega(x,By)=0$ for all $x,y\in IL_0$.
This is equivalent to $(BI)^*=BI$, where $*$ is defined with respect to the
hermitean metric of $V$. 
Thus we can parametrize $\cL_{L_0}$ by the symmetric endomorphisms of $L_0$.
Writing $V=L_0\oplus IL_0$
it is easy to check that 
$$A:=\left(\begin{array}{cc}1&B\\0&1\end{array}\right)\in Sp(V)\ ,$$
$AL_0=L_0$ and $AIL_0=L_1$. Thus $Sp(V)$ acts transitively
on the set $\Lambda^2$ of pairs $(L_0,L_1)$ of transverse Lagrangian
subspaces. 

Let $G$ denote the stabilizer of the pair $(L_0,IL_0)$. 
Let $j:Gl(L_0)\rightarrow Gl(IL_0)$ denote the unique
isomorphism such that 
\begin{equation}\label{rt1}Gl(L_0)\ni A\mapsto \left( \begin{array}{cc}A&0\\0&j(A)\end{array}\right)\in G\ .\end{equation} Then $j(A)=-I(A^{-1})^*I$.
If $L_2\in \Lambda$ is transverse to $L_0$ and $IL_0$, then we write
$$L_2=\left(\begin{array}{cc}1&B\\0&1\end{array}\right) IL_0$$
for some invertible $B$ as above. If $g\in G$ is represented by $A\in Gl(L_0)$
according to (\ref{rt1}),
then 
$$gL_2=\left(\begin{array}{cc}1&-ABIA^*I\\0&1\end{array}\right) IL_0\ .$$
The action of $G$ on $\cL_{L_0}$ is hence given by
action of $Gl(L_0)$ on the symmetric endomorphisms of $L_0$ by
conjugation.
Thus the signature of the symmetric $BI$ is the only invariant of the orbit 
of $G$ generated by $L_2$ inside the space of Lagrangian subspaces
which are transverse to $L_0$ and $IL_0$.
We conclude that $\Lambda^3$ is the disjoint union of orbits
of $Sp(V)$ of points $(L_0,IL_0,L_2)$, which are distinguished
by the signature of a matrix $BI$ defined by $L_2$.

We now consider the orbit generated by a triple $(L_0,IL_0,L_2)\in\Lambda^3$.
The stabilizer $U$ of $(L_0,IL_0,L_2)$ can be identified
with the subgroup of $Gl(L_0)$ fixing the hermitean form on $L_0$ defined by
$BI$. Let $K\subset U$ denote a maximal compact subgroup.
We can choose $K$ such that it fixes the metric $(.,.)$,
hence $K$ is a subgroup of the unitary group of $V$. 
But then it fixes $I$, too. 

Using the explicit formulas given above
one can check that the definition of the closed path  
$\hat{\gamma}(L_0,L_1,L_2)$ only depends on $I$ (and not on 
$\sigma_s,\sigma_t$ or on the choice of coordinates).
It follows that $\hat{\gamma}(L_0,IL_0,L_2)=k\hat{\gamma}(L_0,IL_0,L_2)$ for all $k\in K$.
We now globally trivialize $S,E$ using the parallel transport given by the globally flat metric of $\hat{M}$.
Along the boundary $\partial \hat{M}$ the old and the new trivialization
of $S$ are related by a twist of $-2\pi$ in the structure group of $S$. 
Note that $S^\pm$ are the $\pm\imath$ eigenspaces
of the Clifford multiplication by the volume form of $\hat{M}$.
The image of the path $\hat{\gamma}(L_0,L_1,L_2)$ in the new trivialization
can be obtained (up to homotopy) by $\gamma(L_0,L_1,L_2)(z):=z^{-\imath I} \hat{\gamma}(L_0,L_1,L_2)(z)$,
$z\in S^1=\partial \hat{M}$. We see that $\gamma(L_0,L_1,L_2)$
is $K$-invariant, too.

Let $(x,y)$ be oriented, flat orthonormal coordinates on $\hat{M}$
and write $\hat{\dirac}=\sigma_x \partial_x +\sigma_y \partial_y$.
If we let $K$ act on $E^-$ by, say, $K\ni k\mapsto -\sigma_x^+k\sigma_x^-\in \End(E^-)$, then $\hat{\dirac^+}$ is $K$-equivariant.
We now consider the family of Dirac operators parametrized
by $Sp(V)$, given by
$\hat{\dirac^+}$ subject to the boundary conditions
$Sp(V)\ni g\mapsto \gamma(gL_0,gIL_0,gL_2)$. This family is $K$-equivariant
and we go over to the quotient family
parametrized by $Sp(V)/K$ which we denote by $\check{\dirac}^+$.

Let $X:=\ind(\hat{\dirac}^+)\in R(K)$ be the 
$K$-equivariant index of $\hat{\dirac}^+$ subject to the
boundary conditions given by $\gamma(L_0,IL_0,L_2)$, where $R(K)$ denotes the
representation ring of $K$.
Then $\ind(\check{\dirac}^+)=[Sp(V)\times_K X]\in K^0(Sp(V)/K)$.
The following Lemma implies the proposition for the family $\check{\dirac}^+$.
\begin{lem}\label{ggg}
Let $(n,m)$ be the signature of the quadratic form defined by $BI$ on $L_0$.
Then $X=\C^n-\C^m$, where $\C^n,\C^m$ are the $m$- and $n$-dimensional standard
representations of the corresponding factors of $K\cong U(n)\times U(m)$.
\end{lem}
\proof
Let $T\subset K$ be a maximal torus. If $Y\in R(K)$, then $Y_T$ denotes the restriction to $T$. It is sufficient to show that $X_T=\C^n_T-\C^m_T$.

We first consider the case that $W\cong \C$.
Then $\hat{\dirac}^+$ can be expressed in terms of complex geometry.
Indeed we have 
$$\hat{\dirac}^+=\left(\begin{array}{cc}-\partial&0\\0&\bar{\partial}\end{array}\right)\ .$$ 
Writing $\bar{\partial}=\partial_x + \imath \partial_y$,
$\partial = \partial_x - \imath\partial_y$  we obtain
$$\sigma_x^+=\left(\begin{array}{cc}-1&0\\0&1\end{array}\right),\quad \sigma_y^+=\left(\begin{array}{cc}\imath&0\\0&\imath\end{array}\right)$$
and hence
$$I=\left(\begin{array}{cc}\imath&0\\0&-\imath\end{array}\right)\ .$$
A one-dimensional subspace $\C(a,b)\subset V$, $a,b\in\C$, is Lagrangian
iff $|a|=|b|$. We parametrize the Lagrangian subspaces
of $V$ by $S^1$ associating to $\phi\in S^1$ the space $\C(\phi,1)$.

The space of pairwise transverse triples $\Lambda^3$ consists of
two components which are distinguished by the cyclic order
of the parameters $\phi_i\in S^1$ of $L_i$.
Let $\Lambda^3_1$ be the component with order $\phi_0<\phi_1<\phi_2$ and $\Lambda^3_{-1}$ be the component with order $\phi_0<\phi_2<\phi_1$.
The $T$-equivariant index of $\hat{\dirac^+}$ only depends on the homotopy class of the path $\gamma(L_0,L_1,L_2)$. 

We fix an identification $\pi_1(S^1)=\Z$, $\gamma\mapsto [\gamma]$, such that the path mapped to
$1$ has positive orientation.
We leave to the reader to compute  $[\gamma(L_0,L_1,L_2)]\in \Z$ for
$(L_0,L_1,L_2)\in\Lambda^3_{\pm 1}$.
The result is $[\gamma(L_0,L_1,L_2)]=0$ on $\Lambda^3_1$
and $[\gamma(L_0,L_1,L_2)]=2$ on $\Lambda^3_{-1}$.
\begin{lem}\label{compus}
Let $W\cong\C$, and let the boundary condition of $\hat{\dirac^+}$
be given by a closed path $\gamma$ of Lagrangian subspaces.
Then $\ind(\hat{\dirac^+})=-[\gamma]+1$.
\end{lem}
\proof
Let $n\in\Z$ be represented by the path 
$\gamma_n(\phi)=\C(\phi^n,1)$, $\phi\in S^1$.
The kernel of $\hat{\dirac^+}$ with boundary condition
given by $\gamma_n$ can be identified with the space of
pairs $(f,g)$ of functions on $\hat{M}$, where $\bar{\partial}g=0$, $\partial f=0$
and $z^n g(z)=f(z)$ at $S^1$. This space is non-trivial for $n\le 0$, and it is spanned
by $(\bar{z}^{-n},1), (\bar{z}^{-n-1},z),\dots,(1,z^{-n})$.
One can check that $(\hat{\dirac}^+)^*$ is given by
$$(\dirac^+)^*=\left(\begin{array}{cc}\bar{\partial}&0\\0&-\partial\end{array}\right)\ ,$$ and that the boundary condition is given by the path
$\gamma_{n+2}$. The kernel of $(\hat{\dirac}^+)^*$ can be identified with the
space of pairs $(f,g)$ of functions on $\hat{M}$ with $\bar{\partial}f=0$, $\partial g=0$,
and $z^{n-2}g(z)=f(z)$. The kernel of $(\hat{\dirac}^+)^*$ is non-trivial for 
$n\ge 2$, and it is spanned by
by $(z^{n-2},1),(z^{n-3},\bar{z}),\dots,(1,\bar{z}^{n-2})$.
It follows that $\ind(\hat{\dirac^+})=-n+1$. \hB

We now finish the proof of Lemma \ref{ggg} in the case $W=\C$.
We must show that the signature of the quadratic form given by $BI$ is
$(1,0)$ on $\Lambda^3_1$ and $(0,1)$ on $\Lambda^3_{-1}$.
Let $l_0=(1,1)$, $l_1=(-1,1)$, and $l_2=(-i,1)$
generate the Lagrangian subspaces $L_0,IL_0,L_2$.
Then $(L_0,IL_0,L_2)\in\Lambda^3_1$.
We have $l_2\sim l_0 + \frac{1+\imath}{1-\imath}l_1$ and $Il_0=-il_1$.
It follows that $BI=1$.
The other case is similar. 

In order to complete the proof of Lemma \ref{ggg} in the general case
one reduces to the special case $W=\C$ by considering the direct sums. \hB

We now finish the proof of the proposition.
Let $\{\Lambda^3_i\}$ denote the set of components of $\Lambda^3$ and
choose $x_i\in \Lambda^3_i$ for all $i$. Let $U_i$ be the
stabilizer of $x_i$ in $Sp(V)$.
Then $\bigcup_i Sp(V)/U_i\cong \Lambda^3$ parametrizes the universal family 
of boundary conditions for $\hat{\dirac}^+$ 
given by  a family of path' $B\ni b\mapsto \gamma(L_0(b),L_1(b),L_2(b))$.
Indeed, any such family can be pulled back
from the universal one using the canonical map $B\rightarrow \bigcup_i Sp(V)/U_i$.
Since the fibres of $\pi_i:Sp(V)/K\rightarrow Sp(V)/U_i$ are
symmetric spaces of non-compact type, $\pi_i$ is a homotopy equivalence.
Since the proposition is proved for $\check{\dirac}^+$ (for each component separately) it is also true for the universal family, and hence in general.\hB

\newcommand{\diag}{{\rm diag}}
\newcommand{\arga}{{\rm arg}}
\underline{Remark :}
We sketch another proof of Proposition \ref{expl}
which avoids the use of the universal family and the
equivariant index.

We first consider a model case where $V_{model}:=\C^2$ and $I:=\diag(\imath,-\imath)$.
Let $l_i$, $i=0,1,2$, be the Lagrangian subspaces of $V_{model}$ parametrized by
$\phi_0:=0$, $\phi_1:=\pm\arga(\frac{1-i}{1+i})$, $\phi_2:=-\phi_1$, i.e.,
$l_i=\{(\ee^{\imath \phi_i}x,x)\in\C^2\:|\:x\in\C\}$.
Let $\dirac_{\pm}^+$ be the Dirac operator on $M$ with boundary
conditions given by the triple $l_0,l_1,l_2$.
We have $\ind(\dirac_{\pm}^+)=\pm 1$. This can be proved
in the same way as Proposition \ref{expl}, but the
proof simplifies due to the facts that $\dim(V_{model})=2$, and that
$B$ is a point. 

We now turn to the general case. Let $V$ be any finite dimensional
Hilbert space with hermitean symplectic structure, 
and  let $\{L_0,L_1,L_2\}_{b\in B}$ be a family of pairwise transverse Lagrangian subspaces of $V$.
Then we can write $L_i=\{x+A_iIx\:|\:x\in IL_0\}$, $i=1,2$, where $A_i$ are smooth symmetric bundle endomorphisms of the subbundle $L_0\subset B\times V$
such that $A_1-A_2$ is invertible.

We show that this family of triples is homotopic to a family in some standard form.
Consider the family $A(t)_i:=\frac{1}{2}((1-t)(A_i-A_{3-i})+tA_i)$, $t\in [0,1]$. Then $A(t)_1-A(t)_2=A_1-A_2$
is invertible for all $t$. 
Moreover, $A(1)_i=A_i$ and $A(0)_1=A_1-A_2=-A(0)_2$.
Thus up to homotopy we can assume that $A_1=-A_2=A$, where $A$ is invertible.
There is a further index bundle preserving 
homotopy of $A$ to $A/|A|$. Thus we can assume that $A^2=1$.

Let $\dirac^+$ be the family of operators on $M$ defined by the family
$\{L_0,L_1,L_2\}_{b\in B}$ associated with $A$.
Let $L_0^\pm$ be the $\pm 1$-eigenspaces of $A$.
We will define isomorphisms $\Phi:\ker(\dirac^+_+)\otimes L_0^+\oplus \ker(\dirac^+_{-})\otimes L_0^-\cong \ker(\dirac^+)$, $\Psi:\coker(\dirac^+_+)\otimes L_0^+\oplus \coker(\dirac^+_{-})\otimes L_0^-\cong \coker(\dirac^+)$.

First we fix a basis vector $v_{model}\in l_0$.
Then any $v\in L_0^\pm$ defines an unique symplectic embedding $v_*:V_{model}\hookrightarrow V$ such that $v_*(v_{model})= v$. If  $v,w\in L_0^{\pm}$ and $\mu\in \C$, then  we have $(\mu v+w)_*=\mu v_*+w_*$.
One can check that $v_*(\ker(\dirac^+_{\pm}))\subset \ker(\dirac^+)$, 
$v_*(\coker(\dirac^+_{\pm}))\subset \coker(\dirac^+)$.
We define $\Phi(f\otimes v\oplus f^\prime\otimes v^\prime):=v_*(f)+v^\prime_*(f^\prime)$,
$\Psi(g\otimes w\oplus g^\prime\otimes w^\prime):=w_*(g)+w^\prime_*(g^\prime)$.
In follows that
$$\ind(\dirac^+)=\ind(\dirac^+_+)[L_0^+]+\ind(\dirac^+_{-})[L_0^-]= [L_0^+]-[L_0^-]\in K^0(B)\ .$$
This finishes our sketch of an alternative proof of Proposition \ref{expl}.

\underline{Example :} For the purpose of illustration let us consider an example.
Let $V:=\C^4$ equipped with some complex structure $I$ such that
$L_0=\C^2\subset \C^4$ is Lagrangian.
Let $B:=P^2\C$. If $T\rightarrow B$ denotes the
tautological bundle of $B$, then we have an orthogonal splitting of the
trivial bundle $B\times L_0$ as $T\oplus T^\perp$.
For $b\in B$ let $Q_b$ be the quadratic form on $L_0$
given by the matrix ${\rm diag}(1,-1)$ with respect to the splitting
$L_0=T_b\oplus T_b^\perp$. The family of quadratic forms $\{Q_b\}_{b\in B}$
induces a family of Lagrangian subspaces $\{L_2(b)\}_{b\in B}$ such that
$L_2(b)$ is transverse to $L_0,IL_0$ for all $b\in B$.
Thus $\tau(L_0,IL_0,L_2)$ is defined, and we have
$$\tau(L_0,IL_0,L_2)=\ch(T)-\ch(T^*)=2c_1(T)\ .$$
This class is non-trivial.

\underline{Remark :}
It would be desirable to have an explicit formula for the $\eta$-form
generalizing the result of \cite{leschwojciechowski96}.

\bibliographystyle{plain}

 \end{document}